\documentclass[5p,twocolumn,10pt,times]{elsarticle}
\usepackage{amsthm}
\usepackage{graphicx}
\usepackage{animate}
\usepackage{amsmath}
\usepackage{hyperref}
\begin{document}
\baselineskip11pt

\begin{frontmatter}

\title{Branched splines}

\author{Guohui Zhao\corref{cor1}}
\cortext[cor1]{Corresponding author}
\ead{ghzhao@dlut.edu.cn}
\address{School of Mathematical Sciences, Dalian University of Technology, China}

\begin{abstract}
Spline functions have long been used in numerical solution of differential equations. Recently it revives as isogeometric analysis, which offers integration of finite element analysis and NURBS based CAD into a single unified process. Usually many NURBS pieces are needed to build geometrically continuous CAD models. In this paper, we introduce some splines defined on branched covering of sphere, torus or general domains of $R^2$, which are called branched splines in this paper. A single piece of such splines is enough to build some complex CAD models. Multiresolution analysis on surfaces of high genus built from such splines can be carried out naturally. CAD and FEA are integrated directly on such models. A theoretical framework is presented in this paper, together with some simple examples.
\end{abstract}


\end{frontmatter}

%\linenumbers

\section{Introduction}
Tensor-product splines and NURBS are used to interactively design CAD models. These splines are defined on rectangular domains. Usually it needs many pieces of splines to represent a model. Adjacent pieces are defined independently and transit geometrically continuously. For finite analysis, these splines must be converted into meshes, which is very time consuming. Recently isogeometric analysis\cite{book2} is presented to employ complex NURBS geometry in FEA directly. But it is not easy to adapt independently defined and geometrically continuous splines to isogeometric analysis. In this paper, multivariate splines are generalised to be defined on some standard kinds of manifolds, branched covers of sphere, torus or Euclidean planar domains. Since the definition domain of these splines is a manifold, the models built from these splines can be manifolds. Since splines have natural subdivision schemes, multiresolution analysis can be conducted on such spline manifolds\cite{book3,Matthias1995}.

This paper is organised as follws: Section one is the introduction. Section two lists some basic results of multivariate splines. Section three covers branched covers. Section four defines branched splines. Section five gives some examples. A summary follows in section six.

\section{Prior work}

This paper is mainly concerned with an extension of multivariate splines, which can be used as  a generalization of isogeometric analysis. So some basic results of multivariate splines are listed below. For detailed introduction, please refer to the references and the internet. Multivariate spline theory has been studied extensively. It is used in approximation, CAGD, FEM, etc. For the sake of simplicity, we restrict to the bivariate case. Given a polygonal domain $\Omega$, we partition it with irreducible algebraic curves into finite cells. The partition is denoted by $\Delta$. A function $f$ defined on $\Omega$ is a spline of degree $n$ and smooth order $r$ if $f \in c^r(\Omega)$ and on each cell $\Delta_i$ is a polynomial of degree $n$. The set of all splines of degree $n$ and smooth order $r$ is a linear space and is denoted by $s_n^r(\Omega)$.
Below are three basic theorems of multivariate splines, others will be covered in future papers. The interested readers may refer to \cite{book1}.

    Theorem 1   A piecewise polynomial $p(x,y)$ defined on $\Omega$ with respect to $\Delta$ is a spline of degree $n$ and smooth order $r$ if and only if for any adjacent cells $\Delta_i$ and $\Delta_j$ with common partition curve $l_{ij}$
$$ p_{|\Delta_i}-p_{|\Delta_j}=l_{ij}^{r+1}q_{ij}$$
where $q_{ij}$ is a polynomial and is called the smooth cofactor on $l_{ij}$.

    Theorem 2  Around an interior vertex of $\Delta$, smooth cofactors on partitioning curves concurrent at this vertex satisfy the following equation: $\sum l_{ij}^{r+1}q_{ij} \equiv 0$. This is called the conformality equation around this vertex.

   Theorem 3 The solution space of the conformality equation
   $$\sum_{l=1}^{N}(\alpha_lx+\beta_ly)^{r+1}q_{ij}\equiv 0, \alpha_i\beta_j \neq \alpha_j\beta_i \hspace{0.3cm} for \hspace{0.3cm} i\neq j$$
   has dimension
   \begin{eqnarray*}
   \frac{1}{2}(n-r-[(r+1)/(N-1)])_+ \times ((N-1)n-(N+1)r+ &\\
   (N-3)+(N-1)[(r+1)/(N+1)])
   \end{eqnarray*}
   where $[x]$ denotes the largest integer $\leq x$.

\section{Branched covers}
Branched covers have been used in computer graphics\cite{Wei2007,Felix2007,Lawrence2018}. A brief introduction from Wikipedia is listed below.
\subsection{In topology}

Let $X$ be a topological space. A covering space of $X$ is a topological space $C$ together with a continuous surjective map $f: C\to X$,
such that for every $x \in X$, there exists an open neighborhood $U$ of $x$, such that ${\displaystyle f^{-1}(U)}$ is a union of disjoint open sets in $C$, each of which is mapped homeomorphically onto $U$ by $f$.
The map $f$ is called the covering map, $C$ is called a covering space and $X$ the base space of the covering.

Let $X$ and $C$ be two topological spaces. A continuous surjection $f: C\to X$ is called a branched cover if there exists a nowhere dense subset $\Delta \subset X$ such that $f_{f^{-1}(X\ \Delta)} \rightarrow X\setminus \Delta$ is a covering map. The subset $X\setminus \Delta$ is called the regular set of the branched cover, whereas $\Delta$ its singular set.

\subsection{In Riemann surfaces}
A Riemann surface is a connected complex manifold of complex dimension one, or eqivalently, a Riemann surface is an oriented manifold of real dimension two together with a conformal structure. Given a connected Riemann surface $X$, a branched cover of $X$ is a connected Riemann surface $C$ together with a nonconstant holomorphic map $f: C \to X$.

Suppose $f: C\to X$ is a branched cover and $p\in C$. Then there are local charts of $p$ and $f(p)$ such that $f$ is equal to $z\to z^{e_p}$.  If $e_p>1$, we say that $p$ is a ramification point of $f$ with ramification index $e_p$, and $f(p)$ is a branch point of $f$.

\subsection{Riemann-Hurwitz formula and branched triangulation}
For a branched cover $f: C\to X$ of degree $n$ between two compact Riemann surfaces of genus $g(C)$ and $g(X)$, respectivelly, the Riemann-Hurwitz formula holds: $2g(C)-2=N(g(X)-2)+\sum_{p\in C} (e_p-1)$. This formula is proved as follows:

Choose a triangulation of $X$ such that every branch point of $f$ is a vertex of this triangulation. Write $V,E,F$ for the number of vertices, edges and faces of this triangulation. Pull back by $f$ this triangulation to $C$, Write $V^{'},E^{'},F^{'}$ for the number of vertices, edges and faces of the obtained  triangulation on $C$. As is known that
$$ 2g(X)-2=V-E+F,\hspace{0.5cm} and \hspace{0.5cm} 2g(C)-2=V^{'}-E^{'}+F^{'}.$$
Note that $E^{'}=nE$ and $F^{'}=nF$. Above a branch point $q$ the number of preimages on $C$ is $\#(f^{-1}(q))=n-\sum_{f(p)=q} (e_p-1)$. hence $V^{'}=nV-\sum_{p\in C} (e_p-1)$. The Riemann-Hurwitz formula follows from these equalities.

\begin{figure}
  \includegraphics[width=0.8\linewidth]{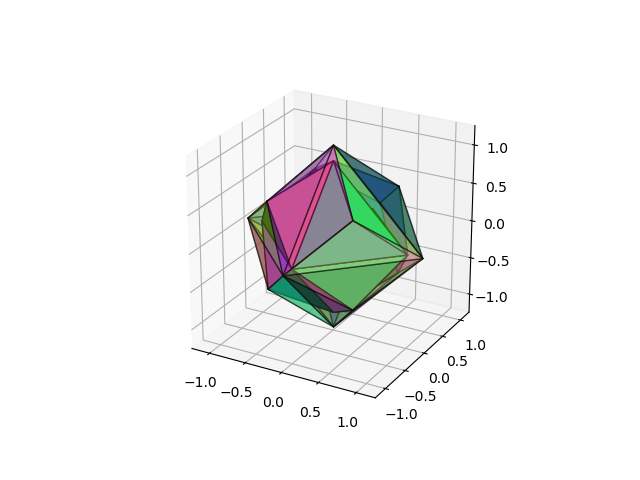}
  \caption{Branched cover}
  \label{fig:cover}
\end{figure}

The pushforward of the triangulation on $C$ in the above proof is called a branched triangulation\ref{fig:cover}, which is a triangulation of a connected manifold formed by $n$ layers of $X$.

%\begin{frame}{Branched cover}
%  \animategraphics[autoplay,loop,controls,width=0.5\linewidth]{10}{genus_3_branched_cover-}{1}{24}
%\end{frame}

\section{Branched splines}
Traditional splines are defined on Euclidean domain, sphere or torus. We generalize them to branched triangulations over Euclidean domain, sphere or torus.

Definition. Given a branched triangulation over Euclidean domain, sphere or torus, we define multivariate splines on the underlying manifold of the branched triangulation just as traditional splines are defined, that is, select a polynomial on each cell of the branched triangulation such that polynomials on adjacent cells satisfy the condition of theorem one. The function obtained in this way is called a branched spline.

If the branched triangulation is over torus, splines on torus, or the traditional biperiodic splines, can be used directly to construct splines on the branched triangulation. The traditional biperiodic B-spline is also a B-spline on the branched triangulation if the interior of its support contains no ramification points. If the interior of the support of the traditional biperiodic B-spline contains one ramification point $p$ or two related ramification points $p$ and $q$, it can be used $e_p$ times resulting in a B-spline on the branched triangulation. Other cases are not clear.

Traditional finite elements can be used on branched triangulation, such as Hsieh-Clough-Tocher triangles,  Powell¨CSabin¨CHeindl triangles and Fraeijs de Veubeke and Sander quadrilaterals.

If the branched triangulation is over a Euclidean domain, branched splines can be constructed similarly and more easily since no periodic conditions are required.

The construction of branched splines on branched triangulations over a sphere is similar, but is postponed to future study since there are no spherical splines at hand, we must construct them first.

\section{Examples}
\begin{figure}
  \includegraphics[width=0.8\linewidth]{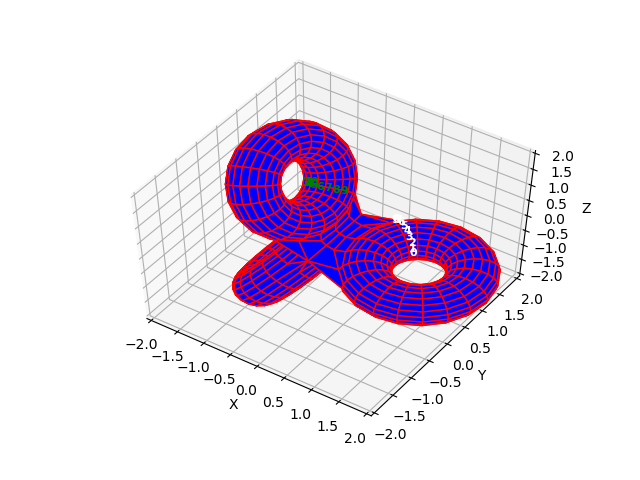}
  \caption{$c^0$ Branched B-spline}
  \label{fig:3-torus1}
\end{figure}

\begin{figure}
  \includegraphics[width=\linewidth]{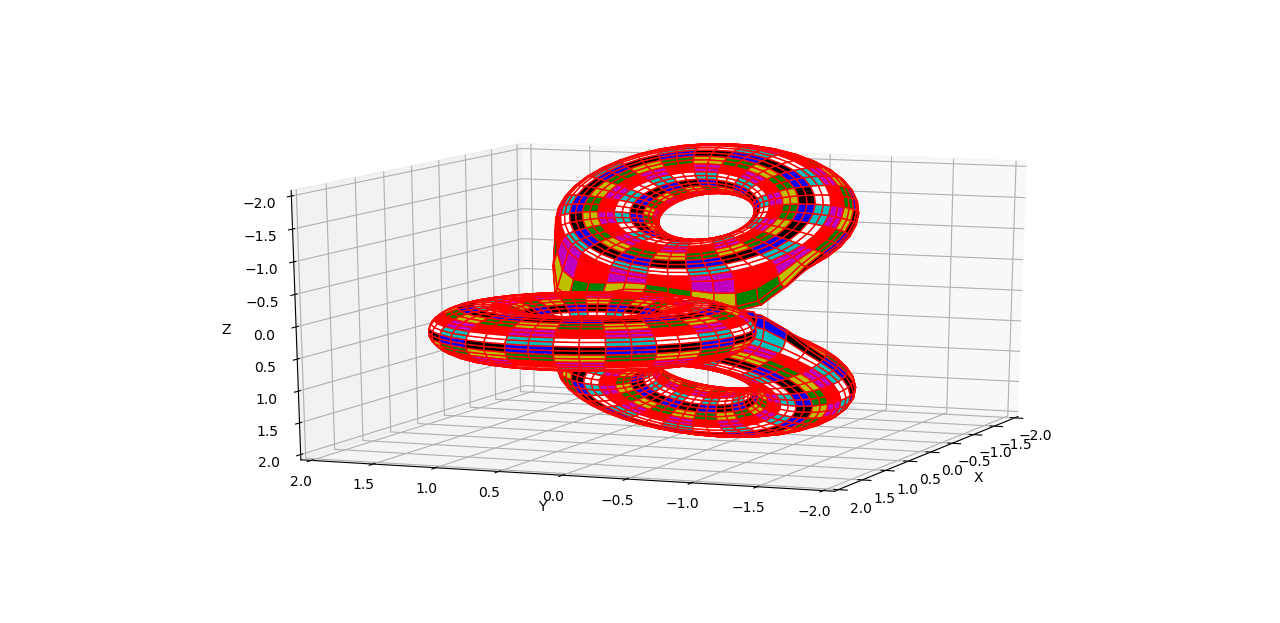}
  \caption{$c^1$ Branched B-spline}
  \label{fig:3-torus2}
\end{figure}

In this section firstly I use branched B-splines on torus to construct a 3-holed torus. The branched triangulation is based on a triple branched cover of torus with two ramification points of ramification index three. The torus is represented by $[0,20]\times[0,20]$ partitioned with lines $x=1,\cdots,19$ and $y=1,\cdots,19$. The ramification points are located at $(10,8)$ and $(10,12)$. $s_1^0$ and $s_2^1$ B-splines are used to construct two examples\ref{fig:3-torus1}\ref{fig:3-torus2}, respectively. The control points are sampled from three tori. The control points of the branched B-splines at the ramification points are the average of three corresponding points on the tori.  The control points of the branched B-splines near the ramification points are modified properly. Interactive design will be introduced later.

\begin{figure}
  \includegraphics[width=\linewidth]{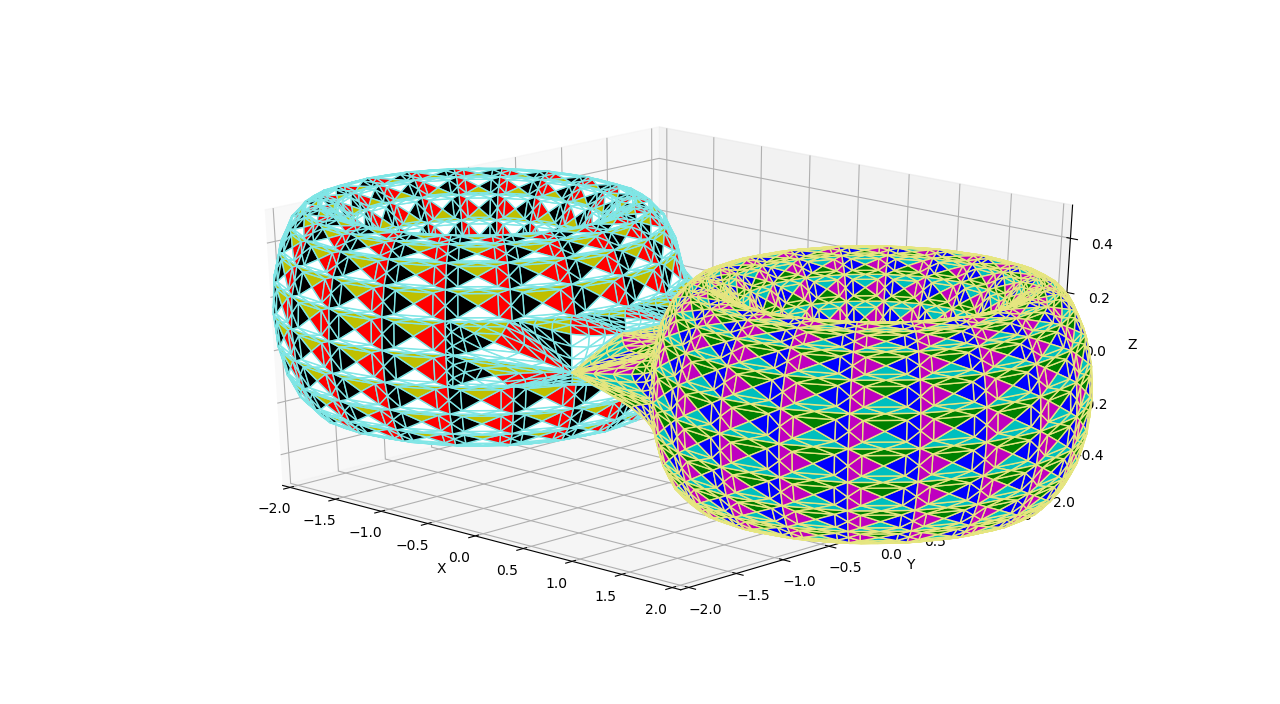}
  \caption{$c^1$ FVS Branched spline}
  \label{fig:4-FVS}
\end{figure}

Secondly I use Fraeijs de Veubeke and Sander quadrilaterals to construct a 2-holed torus. The branched triangulation is based on a double branched cover of torus with two ramification points of ramification index two, which is similar to the one of the above example. An example is given\ref{fig:4-FVS}.

\section*{Summary}
In this paper, branched splines are constructed, a single piece of which can be used to represent higher genus surfaces. So isogeometric analysis can be conducted on such spline surfaces directly. As splines have natural subdivision schemes, multiresolution analysis can be carried out on such spline surfaces, too. These and other numerical analysis problems will be covered in the future.

\end{document}